\documentclass[headsepline,a4paper,openany]{article}
\usepackage{amsthm,amsmath,amssymb,amsfonts}
\usepackage{epsfig}
\usepackage{graphicx,color}
\usepackage{pstricks}
\usepackage{verbatim}
\usepackage[T1]{fontenc}
\usepackage[ansinew]{inputenc}

\newtheorem{lem}{Lemma}
\newtheorem{cor}{Corollary}
\newtheorem{teo}{Theorem}
\newtheorem{cla}{Claim}
\newtheorem{pro}{Proposition}
\newtheorem{obs}{Observation}

\begin{document}

\title{Extremal Graphs for Clique-Paths}
\author{Roman Glebov\\
\small{FU Berlin, Institut f\"{u}r Mathematik, Arnimallee 3, D-14195 Berlin, Germany}\\
\small{glebov@math.fu-berlin.de}}

\maketitle

\begin{abstract}

In this paper we deal with a Tur\'an-type problem: given a
positive integer $n$ and a forbidden graph $H$, how many edges can
there be in a graph on $n$ vertices without a subgraph $H$? 
How does a graph look like if it has this extremal edge number?

The forbidden graph in this article is a clique-path:
a path of length $k$ where each edge is extended to an $r$-clique, $r \geq
3$. We determine both the extremal number and the extremal graphs
for sufficiently large $n$.

\end{abstract}

\section{Introduction}
\label{over}

For integers $n \geq r \geq 1$, we let $T_{n,r}$ denote the \textit{Tur\'an graph}, i.e., the complete
$r$-partite graph on $n$ vertices where each partite set has either
$\left\lfloor n/r \right\rfloor$ or $\left\lceil n/r\right\rceil$ vertices
and the edge set consists of all pairs of vertices from distinct parts. The number
of edges in $T_{n,r}$ is denoted by $t_{n,r}$.
A $K_r$ represents the complete
graph on $r$ vertices.

For a graph $G$ and a vertex $x \in V (G)$, the \textit{neighborhood} of $x$ in $G$ is denoted
by $N_G(x) = \{ y \in V (G) : xy \in E(G)\}$, or if the underlying graph $G$ is clear from the context,
simply $N(x)$.
The neighborhood of a subset $V'$ of vertices is the intersection of the
neighborhoods of the vertices of $V'$, $N(V') = \bigcap_{x \in
V'}N(x)$. The vertices from $N(x)$ are adjacent to $x$, we also
say that $x$ \textit{sees} these vertices.
The \textit{degree} of $x$ in $G$, denoted by $d_G(x)$ or $d(x)$, is the
size of $N_G(x)$. 
We use $\delta(G)$ to denote the minimum
degree in $G$ and $\overline{d}(G)$ for the average degree. 
A vertex $x$ with degree $d(x) =
|V|-1$ is called a \textit{universal vertex}.
If the underlying graph
$G$ is clear from the context, we also write $d_A(x)$ for $|N(x) \cap A|$ with $A \subseteq V (G)$.
For a subset $X \subseteq V (G)$, let $G[X]$ denote the subgraph of $G$ induced by $X$. 
If
$X = V (G)\setminus\{v\}$ for some $v \in V (G)$, we also write $G - v$ for $G[X]$.
A
\textit{matching} in $G$ is a set of edges from $E(G)$, no two of which share a common
vertex.

Suppose that we are given a fixed forbidden graph $H$.
A graph is called \textit{$H$-free}, if it does not contain a copy
of $H$ as a subgraph.
We are interested in the maximum (\textit{extremal}) number,
$ex(n,H)$, of edges an $H$-free graph on $n$ vertices can have.
An $H$-free graph on $n$ vertices with
$ex(n,H)$ edges is called an \textit{extremal graph} for $H$, or just $H$-{\em extremal}.

Mantel~\cite{mantel} determined the extremal number for a triangle, 
and Tur\'an~\cite{turan} generalized
the result and showed that $T_{n,r-1}$ is the unique extremal graph for the
$r$-clique. Although for bipartite graphs even the asymptotics of the extremal
numbers often remains open, Erd{\H{o}}s and Stone~\cite{erdosstone} proved the asymptotical result
$ex(n,H) = (1+o(1))t_{n,\chi(H)-1}$ for non-bipartite graphs $H$. The goal in this
case is now to determine the precise extremal number and all extremal graphs.
Simonovits~\cite{miki} developed a method to find exact extremal numbers using the
stability properties of extremal graphs. A well-known result in this field are for
example the octahedron-free graphs determined by Erd\H{o}s and Simonovits~\cite{octahedron}.

We denote by $F_{k,r}$ the graph on $(r - 1)k + 1$ vertices consisting of $k$ $r$-cliques,
which intersect in exactly one common vertex.

Erd\H{o}s, F\"uredi, Gould and Gunderson \cite{efgg} determined the extremal number
for the $F_{k,3}$ for sufficial large $n$.
Chen, Gould, Pfender and Wei \cite{cgpw} proved the following
generalization of the main theorem of \cite{efgg}:

\begin{teo}
\label{flo} \textit{For every $k \geq 1$ and $r>2$, and for every
$n \geq 16k^3r^8$, if a graph $G$ on $n$ vertices has more than
\[ex(n,K_r) +
\begin{cases}
k^2 - k & \mbox{if $k$ is odd}, \\
k^2 - \frac{3}{ 2}k & \mbox{if $k$ is even}
\end{cases}\]
edges, then $G$ contains a copy of an $F_{k,r}$. Further, the
number of edges is best possible.}
\end{teo}

In this article, we look at $k$ $r$-cliques intersecting in a different
way: let $P_k$ be a $k$-path with $V(P_k)=\{p_1, \dots, p_{k+1}\}$
and $E(P_k)=\{ p_i p_{i+1} : ~ 1\leq i \leq k\}$. 
We extend the edges to $r$-cliques and get a \textit{clique-path} $P_{k,r}$.
Formally, \[V(P_{k,r}):= V(P_k) \cup
\{c_{i,j} : ~ 1\leq i \leq k, 1\leq j \leq r-2\}\] and
\begin{align*}
E(P_{k,r}):= & E(P_k) \cup \{ p_i c_{i,j}: ~ 1\leq i \leq k,
1\leq j \leq r-2\}\\
 & \cup \{ p_{i+1} c_{i,j}: ~ 1\leq i
\leq k, 1\leq j \leq r-2\}\\
& \cup \{ c_{i,a} c_{i,b}: ~ 1\leq i \leq k, 1\leq a \leq b
\leq r-2\} \end{align*}

\begin{figure}[htb]
\begin{center}
\includegraphics[scale= 0.5]{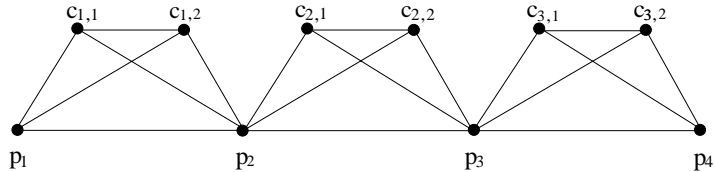}
\caption{A clique-path $P_{3,4}$}
\end{center}
\end{figure}

The aim of this paper is to determine the $P_{k,r}$-extremal graphs.
The graph we want to show being extremal $P_{k,r}$-free on $n$
vertices, $n,k,r$ positive integers with $r \geq 3$ and $n$
sufficiently large, is called $G_{n,k,r}$. It is constructed from
the $(r-1)$-partite Tur\'an graph on
$n-\left\lfloor\frac{k-1}{2}\right\rfloor$ vertices by adding
$\left\lfloor\frac{k-1}{2}\right\rfloor$ universal vertices and
when $k$ is even, also adding an edge. Formally, let
$f:=\left\lfloor\frac{k-1}{2}\right\rfloor$ and let
$G_{n,k,r}:= K_f \vee T_{n-f,r-1}$
be the join of an $f$-clique and the $(r-1)$-partite Tur\'an graph on $n-f$ vertices
with an additional edge, if $k$ is even.

\begin{figure}[htb]
\begin{center}
\includegraphics[scale = 0.55]{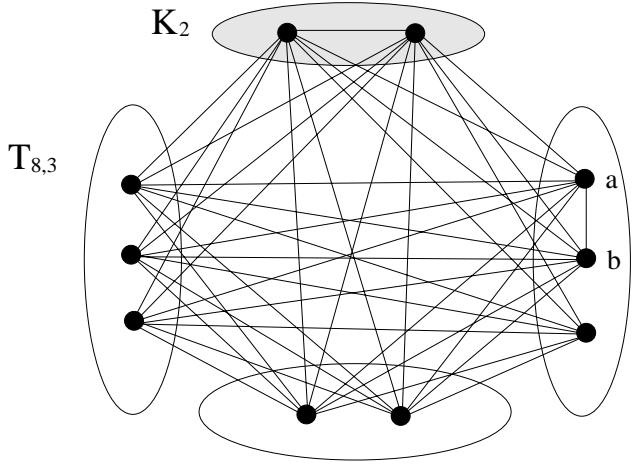}
\caption{A $G_{10,6,4}$}
\end{center}
\end{figure}

Notice that if
$r-1$ is not a factor of $n-f$ and $k$ is even,
there are two nonisomorphic graphs, both called $G_{n,k,r}$,
depending on the size of the set containing $ab$.
Nevertheless, the graphs are quite ``similar'', and the small
difference does not matter in this article, therefore we will not pay
much attention to this fact.

The following is the main theorem of this paper.

\begin{teo}
\label{satz} \textit{Suppose that $G$ is a $P_{k,r}$-free graph on
$n$ vertices with $n > 16k^8r^{11}$.
    Then $|E| \leq g_{n,k,r}$ holds
    and equality occurs if and only if $G$ is isomorphic to a $G_{n,k,r}$.}
\end{teo}

We can easily see that $G_{n,k,r}$ does not contain a $P_{k,r}$.
Since the Tur\'an graph $T_{n-f, r-1}$ is $r-1$-partite, each
$r$-clique in $G_{n,k,r}$ has one vertex from the added $K_f$ or
the edge $ab$. In $P_{k,r}$ there are $k$ $r$-cliques, and each
vertex is contained in two of them at most; the edge $ab$ can
be contained only in one $r$-clique. Thus the length of the
longest clique-path with maximal cliques of order $r$ is
\[k-1 =
\begin{cases}
2\left\lfloor\frac{k-1}{2}\right\rfloor & \mbox{if $k$ is odd}, \\
2\left\lfloor\frac{k-1}{2}\right\rfloor+1 & \mbox{for even $k$.}
\end{cases}\]

The remaining part of the proof of Theorem \ref{satz} is to show that each $P_{k,r}$-free graph on $n$ vertices
with at least $|E(G_{n,k,r})|$ edges is isomorphic to a
$G_{n,k,r}$.

The number of edges in $G_{n,k,r}$ is denoted by
$g_{n,k,r}:= |E(G_{n,k,r})|$.

For the sake of better readability, we omit the graph $G$ in the
notation for the vertex and edge sets and simply write $V$ and $E$
if the underlying graph $G$ is clear from the context. For
simplicity of notation, we will identify isomorphic graphs.

\section{Reduction to high minimum degree}
\label{s2}

The next lemma states the theorem in the case of graphs with high minimum degree.

\begin{lem}
\label{ifdelta} \textit{Suppose that $G$ is a $P_{k,r}$-free graph
on $n$ vertices with $n > 4k^4r^6$,
    and with minimum degree $\delta > \frac{r-2}{r-1}n-1$.
    Then $|E| \leq g_{n,k,r}$ holds
    and equality occurs if and only if $G$ is isomorphic to a $G_{n,k,r}$.}
\end{lem}

We give the proof of this lemma in Section~\ref{main}.
Using the lemma, we prove Theorem~\ref{satz}.

\textit{Proof of Theorem~\ref{satz}.} 
This proof is a so-called ``standard backward induction''.
Let $G$ be a $P_{k,r}$-free graph on $n$ vertices with $n >
16k^8r^{11}$ and $|E| \geq g_{n,k,r}$. Suppose there is a vertex
$x \in V(G)$ with
$d_G(x) < \left\lfloor\frac{r-2}{r-1}n\right\rfloor \leq \delta(G_{n,k,r})$.
We initialize $G^n = G$ and define a process by iteratively deleting vertices with minimum degree.
We continue the process while $\delta(G^i) < \delta(G_{i,k,r})$,
and so during the iterations $|E(G^{i-1})| \geq g_{i-1,k,r} +n-i$.
After $n-l$
steps we get a subgraph $G^l$ with $\delta(G^l) \geq \delta(G_{l,k,r}) \geq \left\lfloor\frac{r-2}{r-1}l\right\rfloor$.
Note that
\[ {l\choose 2} \geq |E(G^l)| \geq g_{l,k,r}
+ n-l > n-l+\frac{r-2}{r-1}\frac{l^2}{2}-\frac{r^2}{2}.\]

Hence $l > \sqrt{rn} \geq 4k^4r^6$ and since $|E(G^{l})| > g_{l,k,r}$,
by Lemma~\ref{ifdelta} $G^l$ contains a $P_{k,r}$. This is a contradiction to the fact that
$G^l$ is $P_{k,r}$-free as a subgraph of a $P_{k,r}$-free graph
$G$.

Thus for all vertices $x \in V(G)$ we obtain $d_G(x) >
\frac{r-2}{r-1}n-1$, hence by Lemma~\ref{ifdelta} $|E| \leq
g_{n,k,r}$ holds and equality occurs if and only if $G$ is
isomorphic to a $G_{n,k,r}$.

\section{The Extremal Graphs for the Clique-Path $P_{k,r}$}
\label{main}

In this section, we prove the remaining claim.

\textit{Proof of Lemma~\ref{ifdelta}.} Suppose that $G$ is a $P_{k,r}$-free graph on $n$ vertices
with $n > 4k^4r^6$, minimum degree $\delta >
\frac{r-2}{r-1}n-1$ and edge number $|E| \geq g_{n,k,r}$. We prove
that $G$ is isomorphic to a $G_{n,k,r}$.

We can assume without loss of generality that $G$ has the most edges
under all graphs that satisfy these properties. We prove the lemma
in a sequence of claims. In the first claim, we see by induction
that the whole graph is close to the Tur\'an graph $T_{n,r-1}$.
In fact it consists of $r-2$ independent sets of size roughly $\frac{n}{r-1}$.
The union of these sets is called $L$. 
In the second claim we show that
excluding only a few vertices from the remaining vertices $V \setminus L$, we
make the edges in that set independent. We call this set $R$ and the
excluded vertices \textit{felons}, since they destroy the structure
of $L$ and $R$. In the third claim and the following proposition
we show some technical statements to prove in Claim~\ref{felons} that there are at most
$\left\lfloor \frac{k-1}{2} \right\rfloor $ of these excluded
felons. The fifth claim says that there is at most one extra edge
inside $L$ or $R$, if $k$ is even, and none for
odd $k$. Then we maximize the number of edges in $G$, and we are
done.

\begin{cla}
\label{l} \textit{$V$ contains a set $L$
that consists of $r-2$ disjoint independent sets, each with $\left\lceil\frac{1}{r-1}n\right\rceil - kr$ vertices.}
\end{cla}

\proof 
In~\cite{turan} and~\cite{cgpw}, the cases
$k=1$ and $k=2$ are already proven. 
We the result for $k=2$ to
start our induction for $k \geq 3$.

\begin{figure}[htb]
\begin{center}
\includegraphics[scale=0.46]{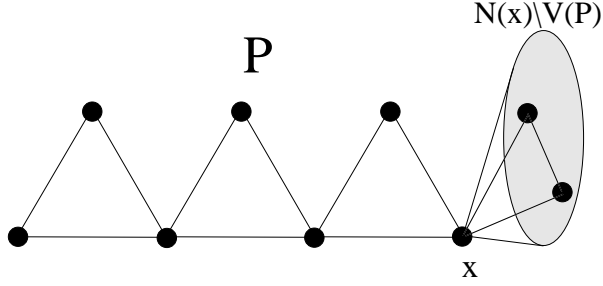}
\caption{The formation of a $P_{4,3}$ from a $P \simeq P_{3,3}$, if $L$ contains a $K_2$.}
\label{3c1}
\end{center}
\end{figure}

Since the function $g_{n,k,r}$ strictly increases with $k$ (either
we get exactly one more edge or one vertex becomes a felon and
gets the full degree $n-1$), by induction, there is a copy $P$ of a $P_{k-1,r}$
as a subgraph in $G$. Let $x$ be the vertex corresponding to $p_k$ in $P$,
and let $L = N(x) \setminus V(P)$ be its neighborhood outside
$P$.
Obviously, $G[L]$ must be $K_{r-1}$-free
since a $K_{r-1}$ in $G[L]$ would extend $P$ into a $P_{k,r}$ via $x$ (see Figure~\ref{3c1}).

Let $l= |L|$ be the
number of vertices in $L$, then
\begin{align}
l &\geq |N(x)| - |V(P)|\nonumber\\
  &> \frac{r-2}{r-1}n-1-(k-1)(r-1)-1 \nonumber\\
  &> \frac{r-2}{r-1}n-kr+r \label{new1}
\end{align}
holds.

Notice that for $r=3$, $L$ is independent, and we are done.
We now assume $r \geq 4$.
We estimate the average degree within $L$:
\begin{align*}
\overline{d}(G[L]) &\geq \delta(G[L]) \\
&\geq \delta(G) -(n-l) \\
&>    (n-\frac{1}{r-1}n-1) -(n-l)  \\
&\overset{\eqref{new1}}{>}    l-\frac{1}{r-3}l \\
&\geq \overline{d}(T_{l,r-3}).
\end{align*}
Thus $|E(G[L])| > ex(l,K_{r-2})$, hence there is an $(r-2)$-clique
$K$ in $G[L]$, $V(K)= \{ v_1, \ldots, v_{r-2}\}$.
We call the $G[L]$-neighborhoods of
$(r-3)$-subsets of $K$
$L_i = N_{G[L]}(K\setminus\{v_i\})$.
Because $G[L]$ is $K_{r-1}$-free, $L_i$'s are independent and pairwise disjoint.
For any vertex $v \in V(K)$, there are less than $\frac{1}{r-1}n+1$
vertices in $L\setminus N(v)$, thus each $L_i$ has more than
\begin{align*}
l-(r-3)\left(\frac{1}{r-1}n+1\right)
&> \frac{1}{r-1}n -kr
\end{align*}
vertices. Together, these neighborhoods form a graph on $r-2$
disjoint independent sets, each with at least $\frac{1}{r-1}n
-kr$ vertices. To establish the claim,
we remove vertices from $L_i$'s
to have exactly $\left\lceil\frac{1}{r-1}n\right\rceil - kr$ vertices in each
and redefine $L := \bigcup L_i$. \qed

Note that
\begin{align}
l = (r-2)\left\lceil\frac{1}{r-1}n\right\rceil - kr(r-2). \label{|l|}
\end{align}

\begin{cla}
\label{r} \textit{$V \setminus L$ can be divided into two parts
$R$ and $F$, $V = R \cup F$, so that the edges in $G[R]$ form a matching
and $|F| \leq k^2r(r-2)+2k<k^2r^2$.}
\end{cla}

\begin{pro}
$G[L]$ contains $k$ pairwise disjoint $(r-2)$-cliques.
\end{pro}

\proof
For $k=3$ note that $l \geq k$, so there are $k$ copies of $K_1$ in $G[L]$.
Now let $k \geq 4$.

Suppose to the contrary that the largest number of pairwise disjoint
$(r-2)$-cliques in $L$ is $k' < k$.
Remove a maximal collection of pairwise disjoint
$(r-2)$-cliques and call the resulting set $L'$.
Then
\begin{align*}
\delta(G[L']) &>  \left (n-\frac{1}{r-1}n-1\right) -(n-l)-kr  \\
    & >   l -\frac{1}{r-3}l  \\
    & \geq \overline{d}(T_{l,r-3}),
\end{align*}
thus $L'$ contains a $K_{r-2}$, leading to a contradiction.
\qed

\textit{Proof of Claim~\ref{r}.}
Note that for each $x \in L_i$,
there are at most $kr$ vertices
outside $L_i$ that are not adjacent to $x$.
Indeed $x$ has at most $\frac{1}{r-1}n +1 < |L_i|+kr$ non-neighbors.

Let us choose a family of $k$ pairwise disjoint $(r-2)$-cliques $C_1, \ldots, C_k$.
Let $F$ be the
set of all vertices in $V \setminus L$,
which are not seen by $\bigcup V(C_i)$,
$F= (V \setminus L)
\setminus N\left(\bigcup V(C_i)\right)$.
We call these vertices \textit{felons}, because they destroy the Tur\'an-like structure
of $G$.
Since for each $v\in L$,  $|(V \setminus L)\setminus N(v)|\leq kr$, we obtain $f:=|F| \leq k^2r(r-2)$. 
Let $R$ be the set of
remaining vertices, $R=V(G) \setminus (L\cup F)$.

Now we remove the vertices of disjoint paths of length at
least $2$ from $R$ greedily until only a matching is left, and add
them to $F$.

\begin{figure}[htb]
\begin{center}
\includegraphics[scale=0.6]{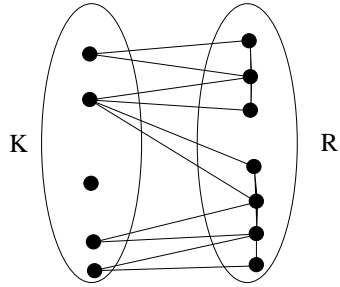}
\caption{The formation of a $P_{5,3}$ from a $P_2$and a $P_3$ in $R$,
where the sum of the lengths of disjoint paths of length at least $2$ in $R$ is at least $5$.}
\label{3c2}
\end{center}
\end{figure}

If the sum of the lengths of all those removed paths was at least
$k$, we could find a $P_{k,r}$ in $G$ similarly to the
Figure~\ref{3c2}.

Thus we removed at most $2k$ vertices. Notice that the edges in
$R$ are pairwise disjoint and there are at most $k^2r(r-2)+2k<k^2r^2$
felons. \qed


\begin{cla}
\label{program} \textit{$V$ can be divided into two parts $S$ and
$F$, where $S$ consists of $r-1$ sets, each with at least
$\frac{1}{r-1}n-k^2r^2$ vertices, where the edges induced by each of them form a matching,
each vertex of $F$ has at least two neighbors in each of these sets,
and $e(S,F) \geq \frac{r-2}{r-1}nf + \left\lfloor\frac{k-1}{2}\right\rfloor \frac{n}{r-1} - \frac{1}{2}\frac{n}{r-1} - k^4
r^4$.}
\end{cla}
\proof
Recall that $L=L_1\cup \ldots \cup L_{r-2}$.
Let $S= R \cup L$ be our set of good vertices. We define $L_{r-1}=R$.
\begin{obs}
\label{|lr|} Notice each of the $L_i$'s has at least
$\frac{1}{r-1}n-k^2r^2+6$ vertices.
\end{obs}
\proof
For $i\leq r-2$, the observation follows from Claim~\ref{l}.
We can see that
\begin{align*}
|R| & \geq n-|L|- k^2r(r-2)-2k \\
& > \frac{1}{r-1}n - k^2r^2 +6.
\end{align*}
\qed

For further calculations, we denote by $e$ the number of edges inside all $G[L_i]$'s,
and fix the current numbers $f' = f$ and $e' = e$, since $f$ and $e$ change soon.
Let us now try to reintegrate some of the felons, that is to insert
them into $S$. To accomplish this, we allow an additional small matching inside the
$L_i$'s and exchange some felons with good vertices. If
there exist a felon $x \in F$ and an $i$ with $d_{L_i}(x) \leq 1$,
add $x$ to $L_i$ and move the resulting degree-2-vertex from $L_i$ to
$F$, should there be one.
Repeat this process, until every felon has at least two neighbors in each $L_i$.

The process terminates in at most $t:=e'+2f'$ steps,
since the value $e+2f$ is reduced by at least $1$
in each iteration.

There are two cases possible for each iteration step:
either one felon and at most one edge are added to an $L_i$,
or we exchange a felon with a good vertex and decrease $e$.
Only in a step where $f$ decreases can $e$ increase by one.
So at most $k^2r(r-2)+2k$ edges are added to the $L_i$'s.
Using Observation~\ref{|lr|}, we have the following observation.
\begin{obs}
\label{claim3} For $i \leq r-2$, still $|L_i| \geq \frac{1}{r-1}n - kr^2$,
and $|L_{r-1}| > \frac{1}{r-1}n - k^2r^2 +6$.
Furthermore, inside each $L_i$ the edges are pairwise disjoint.
\end{obs}

There are at most $\frac{r-2}{r-1}\frac{(n-f)^2}{2}$ edges between
the $L_i$'s,
$\frac{1}{2}\frac{1}{r-1}n+\frac{1}{2}kr^2$ ``old'' edges inside $R$ by (\ref{|l|}) and Observation~\ref{claim3},
$k^2r(r-2)+2k$ edges that came in during the reintegration of felons and
${f \choose 2}$ edges inside $F$. On the other hand,
there are more than $\frac{r-2}{r-1}\frac{n^2}{2} - \frac{r^2}{2} +
\left\lfloor\frac{k-1}{2} \right\rfloor \frac{n}{r-1} -
\frac{k^2}{2}$ edges in $G_{n,k,r}$, thus there are
more than
\begin{align}
&\frac{r-2}{r-1}\frac{n^2}{2} + \left\lfloor \frac{k-1}{2} \right\rfloor
\frac{n}{r-1} - \frac{k^2}{2} - \frac{r^2}{2} \nonumber\\
& \qquad - \frac{r-2}{r-1}\frac{(n-f)^2}{2} - \frac{1}{2}\frac{1}{r-1}n-\frac{1}{2}kr^2 - k^2r(r-2)-2k\nonumber\\
& \geq \frac{r-2}{r-1}nf + \left\lfloor\frac{k-1}{2}\right\rfloor \frac{n}{r-1} - \frac{1}{2}\frac{n}{r-1} - k^4
r^4 \label{eF-S}
\end{align}
edges between $S$ and $F$ in $G$.
\qed

\begin{obs}
\label{unseen} For each $x \in L_i$, there are less than
$k^2r^2$ vertices outside of $L_i$ that are not adjacent to $x$.
\end{obs}

The following two claims will be proven by contradiction using a common principle,
so we state a technical proposition.

For any $P \subset V$ and a clique $K \subset (V\setminus P)\cup F$ with exactly one felon and at most one vertex from
each $L_i$,
an $r$-clique $C$ is a {\em $P$-avoiding extension} of $K$ if
$K \subseteq C$ and $C$ contains (exactly) one vertex from each $L_i \setminus P$.
We use the
following statement:
\begin{pro}
\label{r-clique} For any $P\subset V$ of size $|P| \leq kr$, a $P$-avoiding extension of $K$ exists provided the
felon in $K$ has a degree at least $k^2r^3$ in each
$L_i$ with $L_i \cap K = \emptyset$.
\end{pro}
\proof Let $C$ be a largest clique in $G[(V \setminus P)\cup F]$ containing $K$
with one vertex from each $L_i\setminus P$.
Let $f^{(K)}= K \cap F$ be the felon in $K$.
If there exists an $i$
with $L_i \cap C = \emptyset$, then since $C$ is maximal,
by Observation~\ref{unseen}
\begin{align*}
\left|N_{S}\left(f^{(K)}\right)\right| &\leq \left| \left(L_i \setminus N\left(C\setminus\left\{f^{(K)}\right\}\right)\right) \cup P\right|\\
& < k^2r^3,
\end{align*}
proving the proposition.
\qed

\begin{cla}
\label{felons}
There are at most  $\left\lfloor \frac{k-1}{2}\right\rfloor$ felons,  $f \leq \left\lfloor \frac{k-1}{2}\right\rfloor$.
\end{cla}

\begin{figure}[htb]
\begin{center}
\includegraphics[scale=0.46]{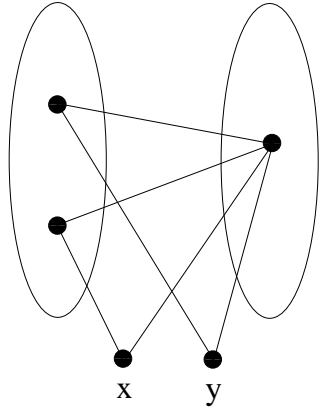}
 \caption{The formation of a connector between two
felons $x$ and $y$ for $r=3$.}
\label{formcon}
\end{center}
\end{figure}

Assume for the sake of contradiction that there are more than
$\left\lfloor \frac{k-1}{2} \right\rfloor$ felons in $G$, $f \geq
\left\lfloor \frac{k+1}{2} \right\rfloor$.
We take the $\left\lfloor \frac{k+1}{2} \right\rfloor$ felons
$f_1, \ldots ,f_{\left\lfloor\frac{k+1}{2} \right\rfloor}$
with the highest $S$-degrees
$d_{S}(f_1) \leq \ldots \leq d_{S}(f_{\left\lfloor\frac{k+1}{2} \right\rfloor})$.
We call $C \simeq P_{2,r}$ a {\em connector} between felons $x$
and $y$, if $x,y \in V(C)$ and $xy \notin E(C)$.
(See Figure~\ref{formcon}.)

\begin{figure}[htb]
\begin{center}
\includegraphics[scale=0.60]{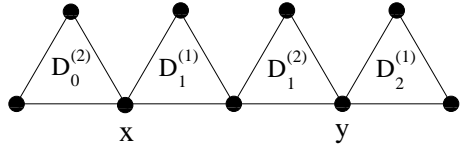}
 \caption{The formation of a $P_{4,3}$ by two
felons $x$ and $y$, their connector and attached triangles.}
\label{conpath}
\end{center}
\end{figure}

Here is our plan: we attach an $r$-clique $D_0^{(2)}$ to $f_1$, then find a
connector with cliques $D_i^{(1)}, D_i^{(2)}$ between $f_i$ and $f_{i+1}$ for $1\leq i \leq
\left\lfloor\frac{k-1}{2}\right\rfloor$.
In the end, we attach an $r$-clique $D_{\left\lfloor\frac{k+1}{2}\right\rfloor}^{(1)}$ to $f_{\left\lfloor\frac{k+1}{2}\right\rfloor}$ and get
the forbidden clique-path $P_{k,r}$,
contradicting the $P_{k,r}$-freeness of $G$. (See Figure~\ref{conpath}.)

This contradiction shall shows that
$|F| \leq \left\lfloor\frac{k-1}{2}\right\rfloor$.

Of course we have to pay attention so the attached connectors and $K_{r-1}$
are only intersecting where they are supposed to.

To be able to use Proposition~\ref{r-clique}, we lower bound the degrees of $f_i$.
By~(\ref{eF-S}) we know that
\[e\left(F \setminus \left\{f_3, \ldots ,f_{\left\lfloor\frac{k+1}{2} \right\rfloor} \right\},S\right)
> \frac{r-2}{r-1}n\left(f-\left\lfloor\frac{k-3}{2}\right\rfloor\right) +
\frac{1}{2}\frac{n}{r-1} - k^4r^4,\]
so since $d_{S}(f_2)$ is maximal within $F \setminus \left\{f_3, \ldots ,f_{\left\lfloor\frac{k+1}{2} \right\rfloor} \right\}$,
we have
\begin{align}
d_{S}(f_2) & > \frac
{\frac{r-2}{r-1}nf + \left\lfloor\frac{k-1}{2}\right\rfloor \frac{n}{r-1} -
\frac{1}{2}\frac{n}{r-1} - k^4 r^4 - n\left\lfloor\frac{k-3}{2}\right\rfloor}
{f-\left\lfloor\frac{k-3}{2}\right\rfloor}\nonumber
\\
& = \frac
{\frac{r-2}{r-1}n\left(f-\left\lfloor\frac{k-3}{2}\right\rfloor\right) +
\frac{1}{2}\frac{n}{r-1} - k^4r^4}{f-\left\lfloor\frac{k-3}{2}\right\rfloor}\nonumber
\\
& > \frac{r-2}{r-1}n + \frac{1}{2}\frac{n}{k^2r^3} - k^2r^2.  \label{dlrf2}
\end{align}
To complete this estimation, the number of neighbors of $f_2$ in
any $L_i$ by Observation~\ref{unseen} is
\begin{align}
d_{L_i}(f_2) & \geq |L_i| - (|S|-d_{S}(f_2)) \nonumber\\
    &\geq \frac{1}{r-1}n - k^2r^2 - \left(n - \frac{r-2}{r-1}n - \frac{1}{2}\frac{n}{k^2r^3} + k^2r^2\right) \nonumber\\
    & = \frac{1}{2}\frac{n}{k^2r^3} - 2k^2r^2\nonumber\\
    & \geq \frac{1}{4}\frac{n}{k^2r^3} - 4k^2r^3 + k^2r^3 \nonumber\\
    & > k^2r^3.  \label{dbf2}
\end{align}
Due to the definition of $f_2$, $d_{L_i}(f_j)> k^2r^3$ for $j \geq 2$.

Since we only need the connectors and the attached cliques to find
the forbidden $P_{k,r}$, we add at most $|V(P_{k,r})|<kr$
vertices to the avoided set $P$. Thus, by~\eqref{dlrf2} and Proposition~\ref{r-clique}
we can make the following corollary.
\begin{cor}
\label{obs} A $P$-avoiding extension always exists for felons $f_i$ with $i \geq
2$.
\end{cor}

Let us choose one neighbor $x$ of $f_1$ in the set $L_j$ with
the fewest neighbors of $f_1$. Because of $\delta_G >
\frac{r-2}{r-1}n-1$, there are less than $\frac{1}{r-1}n+1$
vertices outside of the neighborhood of $f_1$, and at most half of
them (that is $\frac{1}{2}\frac{1}{r-1}n+1$) may be found in each of the
other $L_i$'s with the exception of $L_j$.

Since each of these $L_i$'s has at least
$\frac{1}{r-1}n-k^2r^2+6$ vertices, at least
\begin{align*}
\frac{1}{r-1}n-k^2r^2+6 - \frac{1}{2}\frac{1}{r-1}n-1
& > \frac{1}{2}\frac{1}{r-1}n - k^2r^2 \\
    &\geq  k^2r^3
\end{align*}
of them are in the neighborhood of $f_1$, hence we obtain the following statement.
\begin{obs}
\label{f1}  We can find more than
$k^2r^3$ neighbors of $f_1$ in any $L_i$ with $i \neq j$.
\end{obs}

Hence by Proposition~\ref{r-clique} we can take
a $\emptyset$-avoiding extension $D_0^{(2)}$ of $\{f_1\}$
and call its vertex set $P=V(D_0^{(2)})$.

To construct $D_1^{(1)}$ and $D_1^{(2)}$,
we distinguish two cases:

{\em Case~1}. There is a common neighbor $v$ of $f_1$ and $f_2$ in $L_j\setminus P$.
Then we find a connector between $f_1$ and
$f_2$ the following way.
We take a $P$-avoiding extension $D_1^{(1)}$ of $\{ f_1, v\} $ and
redefine $P := P\cup \left(V\left(D_1^{(1)}\right) \setminus \{v\} \right)$.
Then we take a $P$-avoiding extension $D_1^{(2)}$ of $\{ f_2, v\} $
and add it to $P := P \cup V\left(D_1^{(2)}\right)$.
Using Observation~\ref{f1}
and Proposition~\ref{r-clique}, we find the $D_1^{(1)}$. The existence of
$D_1^{(2)}$ is asserted by Corollary~\ref{obs}.

Obviously, $D_1$ and $D_2$ only intersect in $v$, so they form a
connector between $f_1$ and $f_2$ we searched for.

{\em Case~2}. All the (more than $k^2r^2$) common neighbors of $f_1$ and
$f_2$ in $S$ are outside of $L_j\setminus P$. We proceed similar to the first
case with one difference: we take a neighbor $y$ of $f_1$ in $L_j\setminus P$,
and then find a common neighbor
$v$ of $f_1$ and $f_2$ so that $yv \in E$. To prove that such a vertex exists we make the following proposition.

\begin{pro}
\label{common}
For any $P\subset V$ of size $|P| < kr$ and $i, j \leq \left\lfloor\frac{k+1}{2}\right\rfloor$,
there are more than $k^2r^2$ common neighbors of $f_i$ and
$f_j$ in $S\setminus P$.
\end{pro}
\proof
Using (\ref{dlrf2}) we have
\begin{align}
d_{S\setminus P}(f_i) + d_{S\setminus P}(f_j)
    &\geq d_{S\setminus P}(f_1) + d_{S\setminus P}(f_2) \nonumber\\
	& > d_{S}(f_1) -kr + d_{S}(f_2) - kr \nonumber\\
    &\geq \delta-f + d_{S}(f_2) - 2kr \nonumber\\
    & > \frac{r-2}{r-1}n-1 + \frac{r-2}{r-1}n + \frac{1}{2}\frac{n}{k^2r^3} - 2k^2r^2 -
    2kr \nonumber\\
    & \geq \frac{1}{2}n + \frac{1}{2}n + \frac{1}{2}\frac{n}{k^2r^3} -
    4k^2r^2 + k^2r^2\nonumber\\
    & \geq n + k^2r^2 \nonumber\\
    & > |S| + k^2r^2. \nonumber
\end{align}
\qed

Since by
Observation~\ref{unseen}, there are at most $k^2r^2$ vertices
outside of $L_j$ that are not seen by $y$, we find $v$.
Then we choose a $P$-avoiding extension $D_1^{(1)}$ of the triangle $\{ f_1, y, v\} $,
redefine $P := P\cup \left(V\left(D_1^{(1)}\right) \setminus \{v\} \right)$,
then take a $P$-avoiding extension $D_1^{(2)}$ of $\{ f_2, v\} $,
add it to $P := P \cup v\left(D_1^{(2)}\right)$ and we are done.

It is easier to find connectors between $f_i$ and $f_{i+1}$, $ 2 \leq i \leq
\left\lfloor\frac{k-1}{2}\right\rfloor$.
By Proposition~\ref{common}, we can find a common neighbor $v$ of $f_i$ and
$f_{i+1}$ in $S$, even after having added all the
previous connectors and the attached $r$-clique to $P$ earlier.
Then we take a $P$-avoiding extension $D_i^{(1)}$ of $\{ f_i, v\} $,
redefine $P := P\cup V\left(D_i^{(1)}-v \right)$,
and a $P$-avoiding extension $D_i^{(2)}$ of $\{ f_{i+1}, v\} $,
$P := P \cup V\left(D_i^{(2)}\right)$ and
get a connector between $f_i$ and $f_{i+1}$ following the same argument
as in Case~1.

Finally, for an even $k$ we can attach
a $P$-avoiding extension $D_{\left\lfloor\frac{k+1}{2}\right\rfloor}^{(1)}$
of $\left\{ f_{\left\lfloor\frac{k+1}{2}\right\rfloor}\right\} $ to
$f_{\left\lfloor\frac{k+1}{2}\right\rfloor}$
by Corollary~\ref{obs}.

Hereby we get the forbidden $P_{k,r}$, hence the supposition that
there are at least $\left\lfloor\frac{k+1}{2}\right\rfloor$ felons
was wrong; consequently there are at most
$\left\lfloor\frac{k-1}{2}\right\rfloor$ of them. \qed

\begin{cla}
\label{lr} \textit{There is at most one edge inside the $L_i$'s, if
$k$ is even, and the $L_i$'s are independent for odd $k$.}
\end{cla}

\proof By Claim~\ref{program}, there are at least $\frac{r-2}{r-1}nf +
\left\lfloor\frac{k-1}{2}\right\rfloor \frac{n}{r-1} -
\frac{1}{2}\frac{n}{r-1} - k^4 r^4$ edges between $F$ and
$S$, thus there are at least
$\left\lfloor\frac{k-1}{2}\right\rfloor$ felons. Hence, by
Claim~\ref{program} we have exactly
$\left\lfloor\frac{k-1}{2}\right\rfloor$ felons $f_1, \dots,
f_{\left\lfloor\frac{k-1}{2}\right\rfloor}$.
Let us assume $d_{S}(f_1)\leq \dots \leq d_{S}\left(f_{\left\lfloor\frac{k-1}{2}\right\rfloor}\right)$.
Now $f_1$ has at least
an $S$-degree $d_{S}(f_1) > \frac{r-2}{r-1}n +
\frac{1}{2}\frac{n}{r-1} - k^4 r^4$ (and so do all the other $f_i$'s).
Thus, any two felons $f_i$ and $f_j$ have more than $kr+4$ common neighbors.
Hence, we may delete any four vertices from $S$, initialize $P = \emptyset$
and do a construction similar to the one in Claim~\ref{felons}.
For $i$ from $1$ to $\left\lfloor\frac{k-3}{2}\right\rfloor$, we take
a common neighbor $v$ of $f_i$ and $f_{i+1}$ in $S\setminus P$,
find a $P$-avoiding extension $D_i^{(1)}$ of $\{ f_i, v\} $,
redefine $P := P\cup V\left(D_i^{(1)}-v\right )$,
and then take a $P$-avoiding extension $D_i^{(2)}$ of $\{ f_{i+1}, v\} $
and add it to $P := P \cup V\left(D_i^{(2)}\right)$.
Again the existence of the extensions is asserted by Proposition~\ref{r-clique}.
This way, we constructed a $P$ with
$P_{2\left\lfloor\frac{k-1}{2}\right\rfloor-2, r}\subseteq G[P]$ only with
connectors from each $f_i$ to the successor $f_{i+1}$.

\begin{figure}[htb]
\begin{center}
\includegraphics[scale=0.46]{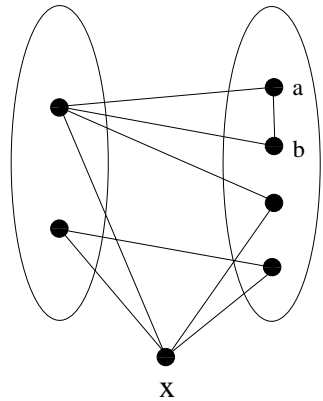}
\caption{The formation of a $P_{3,3}$ with a (red) felon and a (green) extra edge.}
\label{edgepth}
\end{center}
\end{figure}

Our aim is to show that there are not too many edges inside the $L_i$'s. We
only deal with odd $k$ now, since the case of an even $k$
is similar. So what happens if we would find an edge $ab$
inside an $L_i$? We can deport $a$ to the felons, $L_i:= L_i
\setminus \{ a\}$, $F:= F \cup \{ a\}$, and find a common
neighbor $v$ of $a$, $b$ and $f_1$ in $S$. 
Now we attach a connector to our clique-path by
finding a $P$-avoiding extension $D_{0}^{(1)}$
of $\{ f_{1}, v\} $, redefining $P := P \cup V\left(D_{0}^{(1)}- v\right)$,
taking a $P$-avoiding extension $D_{0}^{(2)}$
of $\{a, b,v\}$, and modifying $P := P \cup V\left(D_{0}^{(2)}\right)$.
(See Figure~\ref{edgepth}.)

Now we can attach
a $P$-avoiding extension $D_{\left\lfloor\frac{k-1}{2}\right\rfloor}^{(1)}$
of $\{f_{\left\lfloor\frac{k-1}{2}\right\rfloor}\}$ to the clique-path
and get a forbidden $P_{k,r}$.
But this, as we have already seen, leads us to a contradiction.
Thus the $L_i$s are independent.

Similarly to
that we get at most one edge in the $L_i$s for an even $k$.
\qed

Since $G$ has as many edges as possible, there have to be all the
edges inside $F$ and between $F$ and $S$, and
$S$ has to be a $T_{n-f, r-1}$ (and, of course, the one edge must be
there if $k$ is even), thus $G$ is isomorphic to a $G_{n,k,r}$.
 \qed

To avoid tedious calculations, I did not attempt to lower the
bound $n \geq 16k^8r^{11}$ in the proof, although I strongly
believe the bound can be lowered substantially.

\bigskip

{\bf Acknowledgement.}
My thanks are due to Florian Pfender for suggesting this problem and also for several fruitful comments and discussions.
I would also like to thank Konrad Engel, Thomas Kalinowski, and Tibor Szab\'o, whose com-
ments helped improve this paper considerably.

\end{document}